# Simultaneous estimation of the mean and the variance in heteroscedastic Gaussian regression

### Xavier Gendre

*Laboratoire J.-A. Dieudonné, Université de Nice - Sophia Antipolis*
*Parc Valrose*
*06108 Nice Cedex 02*
*France*
*e-mail:* gendre@unice.fr
*url:* http://math1.unice.fr/∼gendre

**Abstract:** Let $Y$ be a Gaussian vector of $\mathbb{R}^n$ of mean $s$ and diagonal covariance matrix $\Gamma$. Our aim is to estimate both $s$ and the entries $\sigma_i = \Gamma_{i,i}$, for $i = 1, \ldots, n$, on the basis of the observation of two independent copies of $Y$. Our approach is free of any prior assumption on $s$ but requires that we know some upper bound $\gamma$ on the ratio $\max_i \sigma_i / \min_i \sigma_i$. For example, the choice $\gamma = 1$ corresponds to the homoscedastic case where the components of $Y$ are assumed to have common (unknown) variance. In the opposite, the choice $\gamma > 1$ corresponds to the heteroscedastic case where the variances of the components of $Y$ are allowed to vary within some range. Our estimation strategy is based on model selection. We consider a family $\{S_m \times \Sigma_m, \ m \in \mathcal{M}\}$ of parameter sets where $S_m$ and $\Sigma_m$ are linear spaces. To each $m \in \mathcal{M}$, we associate a pair of estimators $(\hat{s}_m, \hat{\sigma}_m)$ of $(s, \sigma)$ with values in $S_m \times \Sigma_m$. Then we design a model selection procedure in view of selecting some $\hat{m}$ among $\mathcal{M}$ in such a way that the Kullback risk of $(\hat{s}_{\hat{m}}, \hat{\sigma}_{\hat{m}})$ is as close as possible to the minimum of the Kullback risks among the family of estimators $\{(\hat{s}_m, \hat{\sigma}_m), \ m \in \mathcal{M}\}$. Then we derive uniform rates of convergence for the estimator $(\hat{s}_{\hat{m}}, \hat{\sigma}_{\hat{m}})$ over Hölderian balls. Finally, we carry out a simulation study in order to illustrate the performances of our estimators in practice.

**AMS 2000 subject classifications:** 62G08.
**Keywords and phrases:** Gaussian regression, heteroscedasticity, model selection, Kullback risk, convergence rate.

## 1. Introduction

Let us consider the statistical framework given by the distribution of a Gaussian vector $Y$ with mean $s = (s_1, \ldots, s_n)' \in \mathbb{R}^n$ and diagonal covariance matrix

$$\Gamma_\sigma = \begin{pmatrix} \sigma_1 & 0 & \cdots & 0 \\ 0 & \ddots & \ddots & \vdots \\ \vdots & \ddots & \ddots & 0 \\ 0 & \cdots & 0 & \sigma_n \end{pmatrix}$$

where $\sigma = (\sigma_1, \ldots, \sigma_n)' \in (0, \infty)^n$. The vectors $s$ and $\sigma$ are both assumed to be unknown. Hereafter, for any $t = (t_1, \ldots, t_n)' \in \mathbb{R}^n$ and $\tau = (\tau_1, \ldots, \tau_n)' \in (0, \infty)^n$, we denote by $P_{t,\tau}$ the distribution of a Gaussian vector with mean $t$ and covariance matrix $\Gamma_\tau$ and by $\mathcal{K}(P_{s,\sigma}, P_{t,\tau})$ the *Kullback-Leibler divergence* between $P_{s,\sigma}$ and $P_{t,\tau}$,

$$\mathcal{K}(P_{s,\sigma}, P_{t,\tau}) = \frac{1}{2} \sum_{i=1}^n \frac{(s_i - t_i)^2}{\tau_i} + \phi\left(\frac{\tau_i}{\sigma_i}\right) ,$$

where $\phi(u) = \log u + 1/u - 1$, for $u > 0$. Note that, if the $\sigma_i$'s are known and constant, the Kullback-Leibler divergence becomes the squared $L^2$-norm and, in expectation, corresponds to the quadratic risk.

Let us suppose that we observe two independent copies of $Y$, namely $Y^{[1]} = (Y_1^{[1]}, \ldots, Y_n^{[1]})'$ and $Y^{[2]} = (Y_1^{[2]}, \ldots, Y_n^{[2]})'$. Their coordinates can be expanded as

$$Y_i^{[j]} = s_i + \sqrt{\sigma_i}\varepsilon_i^{[j]}, \ i = 1, \ldots, n \text{ and } j = 1, 2 , \tag{1.1}$$

where $\varepsilon^{[1]} = (\varepsilon_1^{[1]}, \ldots, \varepsilon_n^{[1]})'$ and $\varepsilon^{[2]} = (\varepsilon_1^{[2]}, \ldots, \varepsilon_n^{[2]})'$ are two independent standard Gaussian vectors. We are interested here in the estimation of the two vectors $s$ and $\sigma$. Indeed, their behaviors contain substantial knowledge about the phenomenon represented by the distribution of $Y$. We have particularly in mind the case of a variance that stays approximately constant by periods and that can take several values in the proceeding of the observations. Of course, we want to estimate the mean $s$ but, in this particular case, we are also interested in recovering the periods of constancy and the values taken by the variance $\sigma$. The Kullback-Leibler divergence measures the differences between two distributions $P_{s,\sigma}$ and $P_{t,\tau}$. Thus, it allows us to deal with the two estimation problems at the same time. More generally, the aim of this paper is to estimate the pair $(s, \sigma)$ by model selection on the basis of the observation of $Y^{[1]}$ and $Y^{[2]}$.

For this, we introduce a collection $\mathcal{F} = \{S_m \times \Sigma_m, \ m \in \mathcal{M}\}$ of products of linear subspaces of $\mathbb{R}^n$ indexed by a finite or countable set $\mathcal{M}$. In the sequel, these products will be called *models* and, for any $m \in \mathcal{M}$, we will denote by $D_m$ the dimension of $S_m \times \Sigma_m$. To each $m \in \mathcal{M}$, we will associate a pair of estimators $(\hat{s}_m, \hat{\sigma}_m)$ that is similar to the *maximum likelihood estimator* (MLE). It is well known that, if the $\sigma_i$'s are equal, the estimators of the mean and the variance factor given by maximization of the likelihood are independent. This fact does not remain true if the $\sigma_i$'s are not constant. To recover the independence between the estimators of the mean and the variance, we construct them separately from the two independent copies $Y^{[1]}$ and $Y^{[2]}$. For the estimator $\hat{s}_m$ of $s$, we take the MLE based on $Y^{[1]}$ and for the estimator $\hat{\sigma}_m$ of $\sigma$, we take the MLE based on $Y^{[2]}$. Thus, for each $m \in \mathcal{M}$, we have a pair of independent estimators $(\hat{s}_m, \hat{\sigma}_m) = (\hat{s}_m(Y^{[1]}), \hat{\sigma}_m(Y^{[2]}))$ with values in $S_m \times \Sigma_m$. The *Kullback risk* of $(\hat{s}_m, \hat{\sigma}_m)$ is given by $\mathbb{E}[\mathcal{K}(P_{s,\sigma}, P_{\hat{s}_m, \hat{\sigma}_m})]$ and is of order of the sum of two terms,

$$\inf_{(t,\tau) \in S_m \times \Sigma_m} \mathcal{K}(P_{s,\sigma}, P_{t,\tau}) + D_m . \tag{1.2}$$

The first one, called the *bias term*, represents the capacity of $S_m \times \Sigma_m$ to approximate the true value of $(s, \sigma)$. The second, called the *variance term*, is proportional to the dimension of the model and corresponds to the amount of noise that we have to control. To warrant a small risk, these two terms have to be small simultaneously. Indeed, using the Kullback risk as a quality criterion, a good model is one minimizing (1.2) among $\mathcal{F}$. Clearly, the choice of a such model depends on the pair of the unknown parameters $(s, \sigma)$ and make good models unavailable to us. So, we have to construct a procedure to select an index $\hat{m} = \hat{m}(Y^{[1]}, Y^{[2]}) \in \mathcal{M}$ depending on the data only, such that $\mathbb{E}[\mathcal{K}(P_{s,\sigma}, P_{\hat{s}_{\hat{m}}, \hat{\sigma}_{\hat{m}}})]$ is close to the smaller risk

$$R(s, \sigma, \mathcal{F}) = \inf_{m \in \mathcal{M}} \mathbb{E}[\mathcal{K}(P_{s,\sigma}, P_{\hat{s}_m, \hat{\sigma}_m})] \ .$$

The art of *model selection* is precisely to provide procedure solely based on the observations in that way. The classical way consists in minimizing an empirical penalized criterion stochastically close to the risk. Considering the *likelihood function* with respect to $Y^{[1]}$,

$$\forall t \in \mathbb{R}^n, \tau \in (0, \infty)^n, \ \mathcal{L}(t, \tau) = \frac{1}{2} \sum_{i=1}^n \frac{\left(Y_i^{[1]} - t_i\right)^2}{\tau_i} + \log \tau_i \ ,$$

we choose $\hat{m}$ as the minimizer over $\mathcal{M}$ of the penalized likelihood criterion

$$\text{Crit}(m) = \mathcal{L}(\hat{s}_m, \hat{\sigma}_m) + \text{pen}(m) \tag{1.3}$$

where pen is a *penalty* function mapping $\mathcal{M}$ into $\mathbb{R}_+ = [0, \infty)$. In this work, we give a form for the penalty in such a way to obtain a pair of estimators $(\hat{s}_{\hat{m}}, \hat{\sigma}_{\hat{m}})$ with a Kullback risk close to $R(s, \sigma, \mathcal{F})$.

Our approach is free of any prior assumption on $s$ but requires that we know some upper bound $\gamma \geqslant 1$ on the ratio

$$\sigma^*/\sigma_* \leqslant \gamma$$

where $\sigma^*$ (resp. $\sigma_*$) is the maximum (resp. minimum) of the $\sigma_i$'s. The knowledge of $\gamma$ allows us to deal equivalently with two different cases. First, "$\gamma = 1$" corresponds to the *homoscedastic* case where the components of $Y^{[1]}$ and $Y^{[2]}$ are independent with a common variance (*i.e.* $\sigma_i \equiv \sigma$) which can be unknown. On the other side, "$\gamma > 1$" means that the $\sigma_i$'s can be distinct and are allowed to vary within some range. This uncommonness of the variances of the observations is known as the *heteroscedastic* case. Heteroscedasticity arises in many practical situations in which the assumption that the variances of the data are equal is debatable.

The research field of the model selection has known an important development in the last decades and it is beyond the scope of this paper to make an exhaustive historical review of the domain. The interested reader could find a good introduction to model selection in the first chapters of [17]. The first heuristics in the domain are due to Mallows [16] for the estimation of the mean in homoscedastic Gaussian regression with known variance. In more general Gaussian framework with common known variance, Barron *et al.* [7], Birgé and Massart ([9] and [10]) have designed an adaptive model selection procedure to estimate the mean for quadratic risk. They provide non-asymptotic upper bound for the risk of the selected estimator. For bound of order of the smaller risk among the collection of models, this kind of result is called *oracle inequalities*. Baraud [5] has generalized their results to homoscedastic statistical models with non-Gaussian noise admitting moment of order larger than 2 and a known variance. All these results remain true for common unknown variance if some upper bound on it is supposed to be known. Of course, the bigger is this bound, the worst are the results. Assuming that $\gamma$ is known does not imply the knowledge of a such upper bound.

In the homoscedastic Gaussian framework with unknown variance, Akaike has proposed penalties for estimating the mean for quadratic risk (see [1, 2] and [3]). Replacing the variance by a particular estimator in his penalty term, Baraud [5] has obtained oracle inequalities for more general noise than Gaussian and polynomial collection of models. Recently, Baraud, Giraud and Huet [6] have constructed penalties able to take into account the complexity of the collection of models for estimating the mean with quadratic risk in Gaussian homoscedastic model with unknown variance. They have also proved results for the estimation of the mean and the variance factor with Kullback risk. This problem is close to ours and corresponds to the case "$\gamma = 1$". A motivation for the present work was to extend their results to the heteroscedastic case "$\gamma > 1$" in order to get oracle inequalities by minimization of penalized criterion as (1.3). Assuming that the collection of models is not too large, we obtain inequalities with the same flavor up to a logarithmic factor

$$\mathbb{E}[\mathcal{K}(P_{s,\sigma}, P_{\hat{s}_{\hat{m}}, \hat{\sigma}_{\hat{m}}})]$$
$$\leqslant C \inf_{m \in \mathcal{M}} \left\{ \inf_{(t,\tau) \in S_m \times \Sigma_m} \mathcal{K}(P_{s,\sigma}, P_{t,\tau}) + D_m \log^{1+\epsilon} D_m \right\} + R \qquad (1.4)$$

where $C$ and $R$ are positive constants depending in particular on $\gamma$ and $\epsilon$ is a positive parameter.

A non-asymptotic model selection approach for estimation problem in heteroscedastic Gaussian model was studied in few papers only. In the chapter 6 of [4], Arlot estimates the mean in heteroscedastic regression framework but for bounded data. For polynomial collection of models, he uses resampling penalties to get oracle inequalities for quadratic risk. Recently, Galtchouk and Pergamenshchikov [14] have provided an adaptive nonparametric estimation procedure for the mean in a heteroscedastic Gaussian regression model. They obtain an oracle inequality for the quadratic risk under some regularity assumptions. Closer to our problem, Comte and Rozenholc [12] have estimated the pair $(s, \sigma)$. Their estimation procedure is different from ours and it makes the theoretical results difficultly comparable between us. For instance, they proceed in two steps (one for the mean and one for the variance) and they give risk bounds separately for each parameter in $L_2$-norm while we estimate directly the pair $(s, \sigma)$ for Kullback risk.

As described in [8], one of the main advantages of inequalities such as (1.4) is that they allow us to derive uniform convergence rates for the risk of the selected estimator over many classes of smoothness. Considering a collection of histogram models, we provide convergence rates over Hölderian balls. Indeed, for $\alpha_1, \alpha_2 \in (0,1]$, if $s$ is $\alpha_1$-Hölderian and $\sigma$ is $\alpha_2$-Hölderian, we prove that the risk of $(\hat{s}_{\hat{m}}, \hat{\sigma}_{\hat{m}})$ converges with a rate of order of

$$\left(\frac{n}{\log^{1+\epsilon} n}\right)^{-2\alpha/(2\alpha+1)}$$

where $\alpha = \min\{\alpha_1, \alpha_2\}$ is the worst regularity. To compare this rate, we can think of the homoscedastic case with only one observation of $Y$. Indeed, in this case, the optimal rate of convergence in the minimax sense is $n^{-2\alpha/(2\alpha+1)}$ and, up to a logarithmic loss, our rate is comparable to this one. To our knowledge, our results in non-asymptotic estimation of the mean and the variance in heteroscedastic Gaussian model are new.

The paper is organized as follows. The main results are presented in section 2. In section 3, we carry out a simulation study in order to illustrate the performances of our estimators in practice with the Kullback risk and the quadratic risk. The last sections are devoted to the proofs and to some technical results.

## 2. Main results

In a first time, we introduce the collection of models, the estimators and the procedure. Next, we present the main results whose proofs can be found in the section 4. In the sequel, we consider the framework (1.1) and, for the sake of simplicity, we suppose that there exists an integer $k_n \geqslant 0$ such that $n = 2^{k_n}$.

### 2.1. Model collection and estimators

In order to estimate the mean and the variance, we consider linear subspaces of $\mathbb{R}^n$ constructed as follows. Let $\mathcal{M}$ be a countable or finite set. To each $m \in \mathcal{M}$, we associate a regular partition $p_m$ of $\{1, \ldots, 2^{k_n}\}$ given by the $|p_m| = 2^{k_m}$ consecutive blocks

$$\left\{(i-1)2^{k_n-k_m}+1, \ldots, i2^{k_n-k_m}\right\}, \ i = 1, \ldots, |p_m| \ .$$

For any $I \in p_m$ and any $x \in \mathbb{R}^n$, let us denote by $x|_I$ the vector of $\mathbb{R}^{n/|p_m|}$ with coordinates $(x_i)_{i \in I}$. Then, to each $m \in \mathcal{M}$, we also associate a linear subspace $E_m$ of $\mathbb{R}^{n/|p_m|}$ with dimension $1 \leqslant d_m \leqslant 2^{k_n-k_m}$. This set of pairs $(p_m, E_m)$ allows us to construct a collection of models. Hereafter, we identify each $m \in \mathcal{M}$ to its corresponding pair $(p_m, E_m)$.

For any $m = (p_m, E_m) \in \mathcal{M}$, we introduce the subspace $S_m \subset \mathbb{R}^n$ of the $E_m$-piecewise vectors,

$$S_m = \{x \in \mathbb{R}^n \text{ such that } \forall I \in p_m, \ x|_I \in E_m\} \ ,$$

and the subspace $\Sigma_m \subset \mathbb{R}^n$ of the piecewise constant vectors,

$$\Sigma_m = \left\{ \sum_{I \in p_m} g_I \mathbb{1}_I, \ \forall I \in p_m, \ g_I \in \mathbb{R} \right\}.$$

The dimension of $S_m \times \Sigma_m$ is denoted by $D_m = |p_m|(d_m + 1)$. To estimate the pair $(s, \sigma)$, we only deal with models $S_m \times \Sigma_m$ constructed in a such way. More precisely, we consider a collection of products of linear subspaces

$$\mathcal{F} = \{S_m \times \Sigma_m, \ m \in \mathcal{M}\} \tag{2.1}$$

where $\mathcal{M}$ is a set of pairs $(p_m, E_m)$ as above. In the paper, we will often make the following hypothesis on the collection of models:

**(H$_\theta$)** There exists $\theta > 1$ such that

$$\forall m \in \mathcal{M}, \ n \geqslant \frac{\theta}{\theta - 1}(\gamma + 2)D_m .$$

This hypothesis avoids handling models with dimension too great with respect to the number of observations.

Let $m \in \mathcal{M}$, we denote by $\pi_m$ the orthogonal projection on $S_m$. We estimate $(s, \sigma)$ by the pair of independent estimators $(\hat{s}_m, \hat{\sigma}_m) \in S_m \times \Sigma_m$ given by

$$\hat{s}_m = \pi_m Y^{[1]}$$

and

$$\hat{\sigma}_m = \sum_{I \in p_m} \hat{\sigma}_{m,I} \mathbb{1}_I \text{ where } \forall I \in p_m, \ \hat{\sigma}_{m,I} = \frac{1}{|I|} \sum_{i \in I} \left( Y_i^{[2]} - \left(\pi_m Y^{[2]}\right)_i \right)^2 .$$

Thus, we get a collection of estimators $\{(\hat{s}_m, \hat{\sigma}_m), \ m \in \mathcal{M}\}$.

### 2.2. Risk upper bound

We first study the risk on a single model to understand its order. Take an arbitrary $m \in \mathcal{M}$. We define $(s_m, \sigma_m) \in S_m \times \Sigma_m$ by

$$s_m = \pi_m s$$

and

$$\sigma_m = \sum_{I \in p_m} \sigma_{m,I} \mathbb{1}_I \text{ where } \forall I \in p_m, \ \sigma_{m,I} = \frac{1}{|I|} \sum_{i \in I} (s_i - s_{m,i})^2 + \sigma_i .$$

Easy computations proves that the pair $(s_m, \sigma_m)$ reaches the minimum of the Kullback-Leibler divergence on $S_m \times \Sigma_m$,

$$\inf_{(t,\tau) \in S_m \times \Sigma_m} \mathcal{K}(P_{s,\sigma}, P_{t,\tau}) = \mathcal{K}(P_{s,\sigma}, P_{s_m, \sigma_m})$$

$$= \frac{1}{2} \sum_{I \in p_m} \sum_{i \in I} \log\left(\frac{\sigma_{m,I}}{\sigma_i}\right) . \tag{2.2}$$

The next proposition allows us to compare this quantity with the Kullback risk of $(\hat{s}_m, \hat{\sigma}_m)$.

**Proposition 1.** *Let $m \in \mathcal{M}$, if the hypothesis ($\mathbf{H_\theta}$) is fulfilled, then*

$$\mathcal{K}(P_{s,\sigma}, P_{s_m,\sigma_m}) \vee \frac{D_m}{4\gamma} \leqslant \mathbb{E}\left[\mathcal{K}\left(P_{s,\sigma}, P_{\hat{s}_m,\hat{\sigma}_m}\right)\right] \leqslant \mathcal{K}(P_{s,\sigma}, P_{s_m,\sigma_m}) + \kappa\gamma^2\theta^2 D_m$$

*where $\kappa > 1$ is a constant that can be taken equal to $1 + 2e^{-1}$.*

As announced in (1.2), this result shows that the Kullback risk of the pair $(\hat{s}_m, \hat{\sigma}_m)$ is of order of the sum of a bias term $\mathcal{K}(P_{s,\sigma}, P_{s_m,\sigma_m})$ and a variance term which is proportional to $D_m$. Thus, minimizing the Kullback risk $\mathbb{E}\left[\mathcal{K}\left(P_{s,\sigma}, P_{\hat{s}_m,\hat{\sigma}_m}\right)\right]$ among $m \in \mathcal{M}$ corresponds to finding a model that realizes a trade-off between these two terms.

Let pen be a non negative function on $\mathcal{M}$, $\hat{m} \in \mathcal{M}$ is any minimizer of the penalized criterion

$$\hat{m} \in \underset{m \in \mathcal{M}}{\operatorname{argmin}} \left\{ \mathcal{L}\left(\hat{s}_m, \hat{\sigma}_m\right) + \operatorname{pen}(m) \right\} . \tag{2.3}$$

In the sequel, we denote by $(\tilde{s}, \tilde{\sigma}) = (\hat{s}_{\hat{m}}, \hat{\sigma}_{\hat{m}})$ the selected pair of estimators. It satisfies the following result:

**Theorem 2.** *Under the hypothesis ($\mathbf{H_\theta}$), suppose there exist $A, B > 0$ such that, for any $(k, d) \in \mathbb{N}^2$,*

$$M_{k,d} = \operatorname{Card}\left\{m \in \mathcal{M} \text{ such that } |p_m| = 2^k \text{ and } d_m = d\right\} \leqslant A(1+d)^B \tag{2.4}$$

*where $\mathcal{M}$ is the set defined at the beginning of the section 2.1. Moreover, assume that there exist $\delta, \epsilon > 0$ such that*

$$D_m \leqslant \frac{5\delta\gamma n}{\log^{1+\epsilon} n}, \quad \forall m \in \mathcal{M} . \tag{2.5}$$

*If we take*

$$\forall m \in \mathcal{M}, \ \operatorname{pen}(m) = \left(\gamma\theta + \log^{1+\epsilon} D_m\right) D_m \tag{2.6}$$

*then*

$$\mathbb{E}\left[\mathcal{K}\left(P_{s,\sigma}, P_{\tilde{s},\tilde{\sigma}}\right)\right] \leqslant C \inf_{m \in \mathcal{M}} \left\{\mathcal{K}\left(P_{s,\sigma}, P_{s_m,\sigma_m}\right) + D_m \log^{1+\epsilon} D_m\right\} + R \tag{2.7}$$

*where $R = R(\gamma, \theta, A, B, \epsilon, \delta)$ is a positive constant and $C$ can be taken equal to*

$$C = 2\left(1 + \frac{(\kappa\gamma\theta + 1)\gamma\theta}{\log^{1+\epsilon} 2}\right) .$$

The inequality (2.7) is close to an oracle inequality up to a logarithmic factor. Thus, considering the penalty (2.6) whose order is slightly larger than the dimension of the model, the risk of the estimator provided by the criterion (1.3) is comparable to the minimum among the collection of models $\mathcal{F}$.

### 2.3. Convergence rate

One of the main advantages of an inequality as (2.7) is that it gives uniform convergence rates with respect to many well known classes of smoothness. To illustrate this, we consider the particular case of the regression on a fixed design. For example, in the framework (1.1), we suppose that

$$\forall 1 \leqslant i \leqslant n, \ s_i = s_r(i/n) \text{ and } \sigma_i = \sigma_r(i/n),$$

where $s_r$ and $\sigma_r$ are two unknown functions that map $[0,1]$ to $\mathbb{R}$.

In this section, we handle the normalized Kullback-Leibler divergence

$$\mathcal{K}_n(P_{s,\sigma}, P_{t,\tau}) = \frac{1}{n}\mathcal{K}(P_{s,\sigma}, P_{t,\tau}) ,$$

and, for any $\alpha \in (0,1)$ and any $L > 0$, we denote by $\mathcal{H}_\alpha(L)$ the space of the $\alpha$-Hölderian functions with constant $L$ on $[0,1]$,

$$\mathcal{H}_\alpha(L) = \{f : [0,1] \to \mathbb{R} \ : \ \forall x, y \in [0,1], |f(x) - f(y)| \leqslant L|x-y|^\alpha\} .$$

Moreover, we consider a collection of models $\mathcal{F}^{PC}$ as described in the section 2.1 such that, for any $m \in \mathcal{M}$, $E_m$ is the space of dyadic piecewise constant functions on $d_m$ blocks. More precisely, let $m = (p_m, E_m) \in \mathcal{M}$ and consider the regular dyadic partition $p'_m$ with $|p_m|d_m$ blocks that is a refinement of $p_m$. We define $S_m$ as the space of the piecewise constant functions on $p'_m$,

$$S_m = \left\{ f = \sum_{I \in p'_m} f_I \mathbb{1}_I \text{ such that } \forall I \in p'_m, \ f_I \in \mathbb{R} \right\} ,$$

and $\Sigma_m$ as the space of the piecewise constant functions on $p_m$,

$$\Sigma_m = \left\{ g = \sum_{I \in p_m} g_I \mathbb{1}_I \text{ such that } \forall I \in p_m, \ g_I \in \mathbb{R} \right\} .$$

Then, the collection of models that we consider is

$$\mathcal{F}^{PC} = \{S_m \times \Sigma_m, \ m \in \mathcal{M}\} .$$

Note that this collection satisfies (2.4) with $A = 1$ and $B = 0$. The following result gives a uniform convergence rate for $(\tilde{s}, \tilde{\sigma})$ over Hölderian balls.

**Proposition 3.** *Let $\alpha_1, \alpha_2 \in (0,1]$, $L_1, L_2 > 0$ and assume that $(\boldsymbol{H_\theta})$ is fulfilled. Consider the collection of models $\mathcal{F}^{PC}$ and $\delta, \epsilon > 0$ such that, for any $m \in \mathcal{M}$,*

$$D_m \leqslant \frac{5\delta\gamma n}{\log^{1+\epsilon} n} .$$

*Denoting by $(\tilde{s}, \tilde{\sigma})$ the estimator selected via the penalty (2.6), if $n$ satisfies*

$$n \geqslant \left(\frac{2\sigma_*^2}{L_1^2 \sigma_* + L_2^2}\right)^2 \vee e^{4(1+\epsilon)^2}$$

*then*

$$\sup_{(s_r,\sigma_r)\in\mathcal{H}_{\alpha_1}(L_1)\times\mathcal{H}_{\alpha_2}(L_2)} \mathbb{E}\left[\mathcal{K}_n(P_{s,\sigma},P_{\tilde{s},\tilde{\sigma}})\right] \leqslant C\left(\frac{n}{\log^{1+\epsilon}n}\right)^{-2\alpha/(2\alpha+1)} \quad (2.8)$$

where $\alpha = \min\{\alpha_1,\alpha_2\}$ and $C$ is a constant which depends on $\alpha_1$, $\alpha_2$, $L_1$, $L_2$, $\theta$, $\gamma$, $\sigma_*$, $\delta$ and $\epsilon$.

For the estimation of the mean $s$ in quadratic risk with one observation of $Y$, Galtchouk and Pergamenshchikov [14] have computed the heteroscedastic minimax risk. Under some assumptions on the regularity of $\sigma_r$ and assuming that $s_r \in \mathcal{H}_{\alpha_1}(L_1)$, they show that the order of the optimal rate of convergence in minimax sense is $C_{\alpha_1,\sigma}n^{-2\alpha_1/(2\alpha_1+1)}$. Concerning the estimation of the variance vector $\sigma$ in quadratic risk with one observation of $Y$ and unknown mean, Wang et al. [19] have proved that the order of the minimax rate of convergence for the estimation of $\sigma$ is $C_{\alpha_1,\alpha_2}\max\left\{n^{-4\alpha_1}, n^{-2\alpha_2/(2\alpha_2+1)}\right\}$ once $s_r \in \mathcal{H}_{\alpha_1}(L_1)$ and $\sigma_r \in \mathcal{H}_{\alpha_2}(L_2)$. For $\alpha_1,\alpha_2 \in (0,1]$ the maximum of these two rates is of order $n^{-2\alpha/(2\alpha+1)}$ where $\alpha = \min\{\alpha_1,\alpha_2\}$ is the worst among the regularities of $s_r$ and $\sigma_r$. Up to a logarithmic term, the rate of convergence over Hölderian balls given by our procedure recover this rate for the Kullback risk.

## 3. Simulation study

To illustrate our results, we consider the following pairs of functions $(s_r,\sigma_r)$ defined on $[0,1]$ and, for each one, we precise the true value of $\gamma$:

- M1 ($\gamma = 2$)

$$s_r(x) = \begin{cases} 4 & \text{if } 0 \leqslant x < 1/4 \\ 0 & \text{if } 1/4 \leqslant x < 1/2 \\ 2 & \text{if } 1/2 \leqslant x < 3/4 \\ 1 & \text{if } 3/4 \leqslant x \leqslant 1 \end{cases} \quad \text{and} \quad \sigma_r(x) = \begin{cases} 2 & \text{if } 0 \leqslant x < 1/2 \\ 1 & \text{if } 1/2 \leqslant x \leqslant 1 \end{cases},$$

- M2 ($\gamma = 1$)

$$s_r(x) = 1 + \sin(2\pi x + \pi/3) \quad \text{and} \quad \sigma_r(x) = 1 ,$$

- M3 ($\gamma = 7/3$)

$$s_r(x) = 3x/2 \quad \text{and} \quad \sigma_r(x) = 1/2 + 2\sin(4\pi(x \wedge 1/2)^2)/3 ,$$

- M4 ($\gamma = 2$)

$$s_r(x) = 1 + \sin(4\pi(x \wedge 1/2)) \quad \text{and} \quad \sigma_r(x) = (3 + \sin(2\pi x))/2 .$$

In all this section, we consider the collection of models $\mathcal{F}^{PC}$ and we take $n = 1024$ (*i.e.* $k_n = 10$). Let us first present how our procedure performs on the examples with the true value of $\gamma$ for each simulation, $\epsilon = 10^{-2}$ and $\delta = 3$ in the assumption (2.5) and the penalty (2.6) with $\theta = 2$. The estimators are drawn in plain line and the true functions in dotted line.

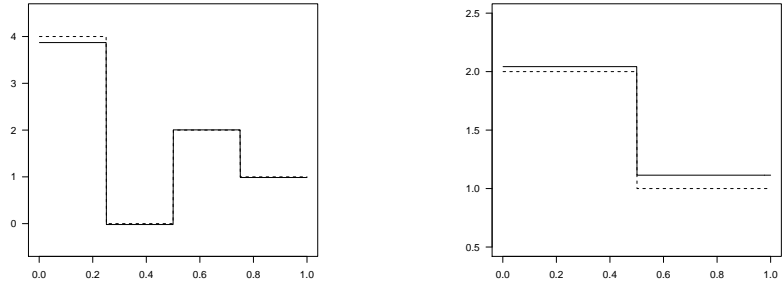

Fig 1: Estimation on the mean (left) and the variance (right) in the case M1.

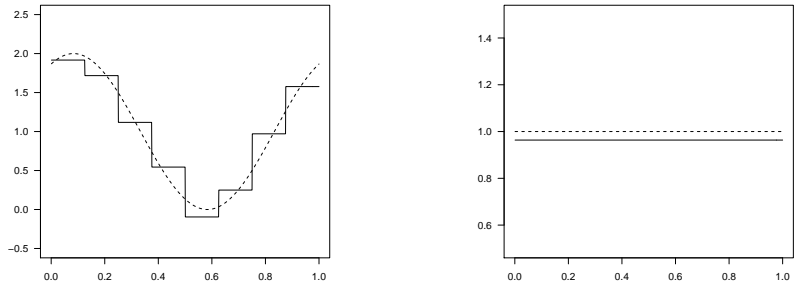

Fig 2: Estimation on the mean (left) and the variance (right) in the case M2.

In the case of M1, we can note that the procedure choose the "good" model in the sense that if the pair $(s_r, \sigma_r)$ belongs to a model of $\mathcal{F}^{PC}$, this one is generally chosen by our procedure. Repeating the simulation 100 000 times with the framework of M1 gives us that, with probability higher than 99.9%, the probability for making this "good" choice is about 0.9978 ($\pm 4 \times 10^{-4}$). Even if the mean does not belong to one of the $S_m$'s, the procedure recover the homoscedastic nature of the observations in the case M2. By doing 100 000 simulations with the framework induced by M2, the probability to choose an homoscedastic model is around 0.99996 ($\pm 1 \times 10^{-5}$) with a confidence of 99.9%. For more general framework as M3 and M4, the estimators perform visually well and detect the changements in the behaviour of the mean and the variance functions.

The parameter $\gamma$ is supposed to be known and is present in the definition of the penalty (2.6). So, we naturally can ask what is its importance in the procedure. In particular, what happens if we do not have the good value? The

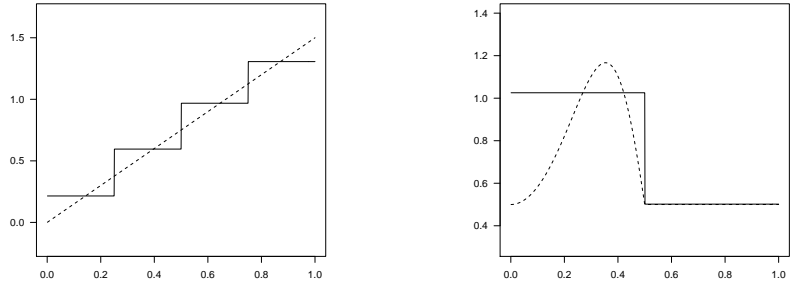

Fig 3: Estimation on the mean (left) and the variance (right) in the case M3.

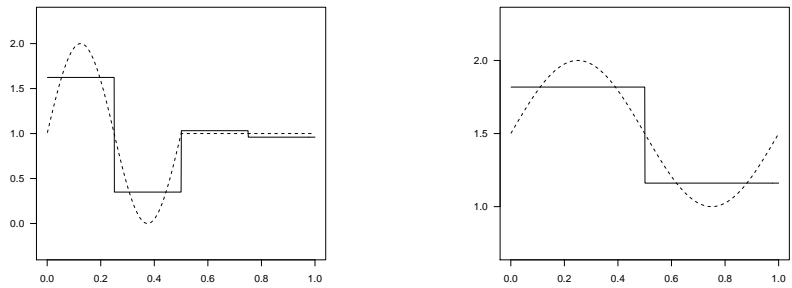

Fig 4: Estimation on the mean (left) and the variance (right) in the case M4.

following table present some estimations of the ratio

$$\mathbb{E}\left[\mathcal{K}\left(P_{s,\sigma}, P_{\tilde{s},\tilde{\sigma}}\right)\right] / \inf_{m \in \mathcal{M}} \mathbb{E}\left[\mathcal{K}\left(P_{s,\sigma}, P_{\hat{s}_m, \hat{\sigma}_m}\right)\right]$$

for several values of $\gamma$. These estimated values have been obtained with 500 repetitions for each one. The main part of the computation time is devoted to the estimation of the oracle's risk. In the cases M1, M3 and M4, the ratio does not suffer to much from small errors on the knowledge of $\gamma$. The more affected case is the homoscedastic one but we see that the best estimation is obtained for the good value of $\gamma$ as we could expect. More generally, it is interesting to observe that, even if there is a small error on the value of $\gamma$, the ratio stays reasonably small.

In the regression framework with heteroscedastic noise, we can be interested in separate estimations of the mean and the variance functions. Because our procedure provide a simultaneous estimation of these two functions, we can ask how perform our estimators $\tilde{s}$ and $\tilde{\sigma}$ individually. Considering the quadratic risks $\mathbb{E}\left[\|s - \tilde{s}\|^2\right]$ and $\mathbb{E}\left[\|\sigma - \tilde{\sigma}\|^2\right]$ of $\tilde{s}$ and $\tilde{\sigma}$ respectively, it could be interesting to compare them to the minimal quadratic risk among the collection of estimators.

| $\gamma$ | 1 | 1.5 | 2 | 2.5 | 3 |
|---|---|---|---|---|---|
| M1 | 0.98 | 1.02 | 1.02 | 1.04 | 1.01 |
| M2 | 1.49 | 1.59 | 1.88 | 2.29 | 2.89 |
| M3 | 1.77 | 1.78 | 1.81 | 1.90 | 1.94 |
| M4 | 1.25 | 1.26 | 1.27 | 1.32 | 1.33 |

Table 1: Ratio between the Kullback risk of $(\tilde{s}, \tilde{\sigma})$ and the one of the oracle

| $\gamma$ | 1 | 1.5 | 2 | 2.5 | 3 |
|---|---|---|---|---|---|
| M1 | 0.98 | 1.01 | 0.95 | 1.04 | 0.98 |
| M2 | 1.52 | 1.67 | 2.04 | 2.43 | 3.04 |
| M3 | 1.73 | 1.76 | 1.82 | 1.88 | 1.96 |
| M4 | 1.47 | 1.48 | 1.47 | 1.47 | 1.49 |

Table 2: Ratio between the $L_2$-risk of $\tilde{s}$ and the minimal one among the $\hat{s}_m$'s

| $\gamma$ | 1 | 1.5 | 2 | 2.5 | 3 |
|---|---|---|---|---|---|
| M1 | 1.00 | 1.06 | 1.03 | 1.02 | 1.01 |
| M2 | 1.11 | 1.56 | 1.68 | 2.21 | 3.36 |
| M3 | 2.02 | 2.07 | 2.13 | 2.20 | 2.23 |
| M4 | 1.18 | 1.37 | 1.34 | 1.44 | 1.49 |

Table 3: Ratio between the $L_2$-risk of $\tilde{\sigma}$ and the minimal one among the $\hat{\sigma}_m$'s

To illustrate this, we give below two sets of estimations of the ratios

$$\mathbb{E}\left[\|s-\tilde{s}\|^2\right] / \inf_{m\in\mathcal{M}} \mathbb{E}\left[\|s-\hat{s}_m\|^2\right] \quad \text{and} \quad \mathbb{E}\left[\|\sigma-\tilde{\sigma}\|^2\right] / \inf_{m\in\mathcal{M}} \mathbb{E}\left[\|\sigma-\hat{\sigma}_m\|^2\right]$$

in the frameworks presented in the beginning of this section. We can observe on the following estimations that the quadratic risks of our estimators are quite close to the minimal ones among the collection of models.

## 4. Proofs

For any $I \subset \{1, \ldots, n\}$ and any $x, y \in \mathbb{R}^n$, we introduce the notations

$$\langle x, y \rangle_I = \sum_{i \in I} x_i y_i \quad \text{and} \quad \|x\|_I^2 = \sum_{i \in I} x_i^2 .$$

Let $m \in \mathcal{M}$, we will use several times in the proofs the fact that, for any $I \in p_m$,

$$|I|\hat{\sigma}_{m,I} \geqslant \sigma_* \chi^2(|I| - d_m - 1) \tag{4.1}$$

where $\chi^2(|I| - d_m - 1)$ is a $\chi^2$ random variable with $|I| - d_m - 1$ degrees of freedom.

### 4.1. Proof of the proposition 1

Recalling (2.2) and using the independence between $\hat{s}_m$ and $\hat{\sigma}_m$, we expand the Kullback risk of $(\hat{s}_m, \hat{\sigma}_m)$,

$$\mathbb{E}[\mathcal{K}(P_{s,\sigma}, P_{\hat{s}_m, \hat{\sigma}_m})] = \frac{1}{2} \sum_{I \in p_m} \sum_{i \in I} \mathbb{E}\left[\frac{(s_i - \hat{s}_{m,i})^2}{\hat{\sigma}_{m,I}} + \phi\left(\frac{\hat{\sigma}_{m,I}}{\sigma_i}\right)\right] \quad (4.2)$$

$$= \frac{1}{2} \sum_{I \in p_m} \mathbb{E}\left[\frac{1}{\hat{\sigma}_{m,I}}\right] \mathbb{E}\left[\|s - \hat{s}_m\|_I^2\right]$$

$$+ \frac{1}{2} \sum_{I \in p_m} \sum_{i \in I} \mathbb{E}\left[\log \frac{\hat{\sigma}_{m,I}}{\sigma_{m,I}} + \frac{\sigma_i}{\hat{\sigma}_{m,I}} - 1\right] + \log \frac{\sigma_{m,I}}{\sigma_i}$$

$$= \mathcal{K}(P_{s,\sigma}, P_{s_m, \sigma_m}) + \frac{1}{2} \sum_{I \in p_m} |I| \mathbb{E}\left[\phi\left(\frac{\hat{\sigma}_{m,I}}{\sigma_{m,I}}\right)\right]$$

$$+ \frac{1}{2} \sum_{I \in p_m} \sum_{i \in I} \mathbb{E}\left[\frac{\sigma_i + (s_i - s_{m,i})^2 - \sigma_{m,I}}{\hat{\sigma}_{m,I}}\right]$$

$$+ \frac{1}{2} \sum_{I \in p_m} \mathbb{E}\left[\frac{1}{\hat{\sigma}_{m,I}}\right] \mathbb{E}\left[\|\pi_m \Gamma_\sigma^{1/2} \varepsilon^{[1]}\|_I^2\right]$$

$$= \mathcal{K}(P_{s,\sigma}, P_{s_m, \sigma_m}) + \mathbb{E}_1 + \mathbb{E}_2 \quad (4.3)$$

where

$$\mathbb{E}_1 = \frac{1}{2} \sum_{I \in p_m} |I| \mathbb{E}\left[\phi\left(\frac{\hat{\sigma}_{m,I}}{\sigma_{m,I}}\right)\right] \quad \text{and} \quad \mathbb{E}_2 = \frac{1}{2} \sum_{I \in p_m} \mathbb{E}\left[\frac{1}{\hat{\sigma}_{m,I}}\right] \sum_{i \in I} \pi_{m,i,i} \sigma_i .$$

To upper bound the first expectation, note that

$$\forall I \in p_m, \ \mathbb{E}[\hat{\sigma}_{m,I}] = \sigma_{m,I} - \frac{1}{|I|} \sum_{i \in I} \pi_{m,i,i} \sigma_i = \sigma_{m,I}(1 - \rho_I)$$

where

$$\rho_I = \frac{1}{|I|\sigma_{m,I}} \sum_{i \in I} \pi_{m,i,i} \sigma_i \in (0, 1) .$$

We apply the lemmas 10 and 11 to each block $I \in p_m$ and, by concavity of the logarithm, we get

$$\mathbb{E}\left[\phi\left(\frac{\hat{\sigma}_{m,I}}{\sigma_{m,I}}\right)\right] \leqslant \log \mathbb{E}\left[\frac{\hat{\sigma}_{m,I}}{\sigma_{m,I}}\right] + \mathbb{E}\left[\frac{\sigma_{m,I}}{\hat{\sigma}_{m,I}}\right] - 1$$

$$\leqslant \log(1 - \rho_I) + \frac{1}{1 - \rho_I}\left(1 + \frac{2\kappa\gamma^2}{|I| - d_m - 2}\right) - 1$$

$$\leqslant -\rho_I + \frac{1}{1 - \rho_I}\left(1 + \frac{2\kappa\gamma^2}{|I| - d_m - 2}\right) - 1$$

$$\leqslant \frac{1}{1 - \rho_I}\left(\rho_I^2 + \frac{2\kappa\gamma^2}{|I| - d_m - 2}\right) .$$

Using $(\mathbf{H}_\theta)$ and the fact that $\rho_I \leqslant \gamma d_m/|I|$, we obtain

$$\begin{aligned}
\mathbb{E}_1 &\leqslant \frac{1}{2} \sum_{I \in p_m} \frac{|I|}{1-\rho_I} \left( \rho_I^2 + \frac{2\kappa\gamma^2}{|I|-d_m-2} \right) \\
&\leqslant \frac{1}{2} \sum_{I \in p_m} \frac{\gamma^2 d_m^2}{|I|-\gamma d_m} + \frac{2\kappa\gamma^2|I|^2}{(|I|-\gamma d_m)(|I|-d_m-2)} \\
&\leqslant \frac{\gamma^2 \theta |p_m| d_m}{2} + \kappa\gamma^2\theta^2 |p_m| \ . \quad (4.4)
\end{aligned}$$

The second expectation in (4.3) is easier to upper bound by using (4.1) and the fact that $d_m \geqslant 1$,

$$\begin{aligned}
\mathbb{E}_2 &= \frac{1}{2} \sum_{I \in p_m} \mathbb{E}\left[\frac{1}{\hat{\sigma}_{m,I}}\right] \sum_{i \in I} \pi_{m,i,i}\sigma_i \\
&\leqslant \frac{1}{2} \sum_{I \in p_m} \frac{\gamma |I| d_m}{|I|-d_m-3} \\
&\leqslant \frac{\gamma \theta |p_m| d_m}{2} \ . \quad (4.5)
\end{aligned}$$

We now sum (4.4) and (4.5) to obtain

$$\mathbb{E}_1 + \mathbb{E}_2 \leqslant \gamma^2 \theta |p_m| d_m + \kappa\gamma^2\theta^2 |p_m| \leqslant \kappa\gamma^2\theta^2 D_m \ .$$

For the lower bound, the positivity of $\phi$ in (4.2) and the independence between $\hat{s}_m$ and $\hat{\sigma}_m$ give us

$$\begin{aligned}
\mathbb{E}\left[\mathcal{K}(P_{s,\sigma}, P_{\hat{s}_m, \hat{\sigma}_m})\right] &\geqslant \frac{1}{2} \sum_{I \in p_m} \mathbb{E}\left[\frac{\|s-\hat{s}_m\|_I^2}{\hat{\sigma}_{m,I}}\right] \\
&\geqslant \frac{1}{2} \sum_{I \in p_m} \frac{\mathbb{E}\left[\|s-\hat{s}_m\|_I^2\right]}{\mathbb{E}\left[\hat{\sigma}_{m,I}\right]} \\
&\geqslant \frac{1}{2} \sum_{I \in p_m} |I| \frac{\|s-s_m\|_I^2 + \sigma_* d_m}{\|s-s_m\|_I^2 + (|I|-d_m)\sigma^*} \ .
\end{aligned}$$

It is obvious that the hypothesis $(\mathbf{H}_\theta)$ ensures $d_m \leqslant |I|/2$. Thus, we get $\sigma_* d_m \leqslant (|I|-d_m)\sigma^*$ and

$$\mathbb{E}\left[\mathcal{K}(P_{s,\sigma}, P_{\hat{s}_m, \hat{\sigma}_m})\right] \geqslant \frac{1}{2} \sum_{I \in p_m} \frac{|I|\sigma_* d_m}{(|I|-d_m)\sigma^*} \geqslant \frac{|p_m| d_m}{2\gamma} \geqslant \frac{D_m}{4\gamma} \ .$$

To conclude, we know that $(\hat{s}_m, \hat{\sigma}_m) \in S_m \times \Sigma_m$ and, by definition of $(s_m, \sigma_m)$, it implies

$$\mathbb{E}\left[\mathcal{K}(P_{s,\sigma}, P_{\hat{s}_m, \hat{\sigma}_m})\right] \geqslant \mathcal{K}(P_{s,\sigma}, P_{s_m, \sigma_m}) \ .$$

### 4.2. Proof of theorem 2

We prove the following more general result:

**Theorem 4.** *Let $\alpha \in (0,1)$ and consider a collection of positive weights $\{x_m\}_{m \in \mathcal{M}}$. If the hypothesis $(\boldsymbol{H_\theta})$ is fulfilled and if*

$$\forall m \in \mathcal{M}, \ pen(m) \geqslant \gamma\theta D_m + x_m \ , \tag{4.6}$$

*then*

$$(1-\alpha)\mathbb{E}\left[\mathcal{K}\left(P_{s,\sigma}, P_{\tilde{s},\tilde{\sigma}}\right)\right]$$
$$\leqslant \inf_{m \in \mathcal{M}}\left\{\mathbb{E}\left[\mathcal{K}\left(P_{s,\sigma}, P_{\hat{s}_m,\hat{\sigma}_m}\right)\right] + pen(m)\right\} + R_1(\mathcal{M}) + R_2(\mathcal{M})$$

*where $R_1(\mathcal{M})$ and $R_2(\mathcal{M})$ are defined by*

$$R_1(\mathcal{M}) = C\theta^2\gamma \sum_{m \in \mathcal{M}} \sqrt{|p_m|}d_m \left(\frac{2C\theta^2\gamma\sqrt{|p_m|}d_m \log(1+d_m)}{x_m}\right)^{\lfloor 2\log(1+d_m)\rfloor}$$

*and*

$$R_2(\mathcal{M}) = \frac{2(\alpha+\gamma\theta)+1}{\alpha} \sum_{m \in \mathcal{M}} |p_m| \exp\left(-\frac{n}{2\theta|p_m|}\log\left(1+\frac{\alpha|p_m|x_m}{\gamma n(\alpha+2)}\right)\right) \ .$$

*In these expressions, $\lfloor \cdot \rfloor$ is the integral part and $C$ is a positive constant that could be taken equal to $12\sqrt{2e}/(\sqrt{e}-1)$.*

Before proving this result, let us see how it implies the theorem 2. The choice (2.6) for the penalty function corresponds to $x_m = D_m \log^{1+\epsilon} D_m$ in (4.6). Applying the previous theorem with $\alpha = 1/2$ leads us to

$$\mathbb{E}\left[\mathcal{K}\left(P_{s,\sigma}, P_{\tilde{s},\tilde{\sigma}}\right)\right]$$
$$\leqslant 2\inf_{m \in \mathcal{M}}\left\{\mathbb{E}\left[\mathcal{K}\left(P_{s,\sigma}, P_{\hat{s}_m,\hat{\sigma}_m}\right)\right] + pen(m)\right\} + 2C\theta^2\gamma R_1 + 8(\gamma\theta+1)R_2$$

with

$$R_1 = \sum_{m \in \mathcal{M}} \sqrt{|p_m|}d_m \left(\frac{2C\theta^2\gamma\sqrt{|p_m|}d_m \log(1+d_m)}{x_m}\right)^{\lfloor 2\log(1+d_m)\rfloor}$$

and

$$R_2 = \sum_{m \in \mathcal{M}} |p_m| \exp\left(-\frac{n}{2\theta|p_m|}\log\left(1+\frac{|p_m|x_m}{5\gamma n}\right)\right) \ .$$

Using the upper bound on the risk of the proposition 1, we easily obtain the coefficient of the infimum in (2.7). Thus, it remains to prove that the two quantities $R_1$ and $R_2$ can be upper bounded independently of $n$. For this, we denote by $B' = B + 2\log(2C\theta^2\gamma) + 1$ and we compute

$$\begin{aligned}
R_1 &= \sum_{m \in \mathcal{M}} \sqrt{|p_m|} d_m \left( \frac{2C\theta^2\gamma\sqrt{|p_m|}d_m \log(1+d_m)}{|p_m|(1+d_m)\log^{1+\epsilon}(|p_m|(1+d_m))} \right)^{\lfloor 2\log(1+d_m) \rfloor} \\
&\leqslant \sum_{k \geqslant 0} \sum_{d \geqslant 1} M_{k,d} 2^{k/2} d \left( 2C\theta^2\gamma 2^{-k/2} \frac{\log(1+d)}{(k\log 2 + \log(1+d))^{1+\epsilon}} \right)^{\lfloor 2\log(1+d) \rfloor} \\
&\leqslant A \sum_{k \geqslant 0} \sum_{d \geqslant 1} (1+d)^{B'} 2^{k/2} \left( \frac{2^{-k/2}\log(1+d)}{(k\log 2 + \log(1+d))^{1+\epsilon}} \right)^{\lfloor 2\log(1+d) \rfloor} \\
&\leqslant A(R'_1 + R''_1) .
\end{aligned}$$

We have split the sum in two terms, the first one is for $d = 1$,

$$R'_1 = \sum_{k \geqslant 0} \frac{2^{B'} \log 2}{(k\log 2 + \log 2)^{1+\epsilon}} = \frac{2^{B'}}{\log^{\epsilon} 2} \sum_{k \geqslant 0} \frac{1}{(k+1)^{1+\epsilon}} < \infty .$$

The other part $R''_1$ is for $d \geqslant 2$ and is equal to

$$\sum_{k \geqslant 0} \sum_{d \geqslant 2} (1+d)^{B'} 2^{-k(\lfloor 2\log(1+d) \rfloor - 1)/2} \left( \frac{\log(1+d)}{(k\log 2 + \log(1+d))^{1+\epsilon}} \right)^{\lfloor 2\log(1+d) \rfloor} .$$

Noting that $1 < \log(1+d) \leqslant \lfloor 2\log(1+d) \rfloor$, we have

$$\begin{aligned}
R''_1 &\leqslant \sum_{k \geqslant 0} 2^{-k/2} \sum_{d \geqslant 2} (1+d)^{B'} \exp\left(-\epsilon \lfloor 2\log(1+d) \rfloor \log\log(1+d)\right) \\
&\leqslant \frac{\sqrt{2}}{\sqrt{2}-1} \sum_{d \geqslant 2} (1+d)^{B' - \epsilon \log\log(1+d)} < \infty .
\end{aligned}$$

We now handle $R_2$. Our choice of $x_m = D_m \log^{1+\epsilon} D_m$ and the hypothesis (2.5) imply

$$\frac{|p_m|x_m}{5\gamma n} \leqslant \delta|p_m| = \frac{1 - (\delta|p_m|+1)^{-1}}{(\delta|p_m|+1)^{-1}} .$$

We recall that, for any $a \in (0,1)$, if $0 \leqslant t \leqslant (1-a)/a$, then $\log(1+t) \geqslant at$. Take $a = (\delta|p_m|+1)^{-1}$ to obtain

$$\log\left(1 + \frac{|p_m|x_m}{5\gamma n}\right) \geqslant \frac{|p_m|x_m}{5(\delta|p_m|+1)\gamma n} \geqslant \frac{x_m}{5(\delta+1)\gamma n} .$$

For any positive $t$, $1 + t^{1+\epsilon} \leqslant (1+t)^{1+\epsilon}$, then we finally obtain

$$
\begin{aligned}
R_2 &= \sum_{m \in \mathcal{M}} |p_m| \exp\left(-\frac{n}{2\theta |p_m|} \log\left(1 + \frac{|p_m| x_m}{5\gamma n}\right)\right) \\
&\leqslant \sum_{m \in \mathcal{M}} |p_m| \exp\left(-\frac{x_m}{10\theta\gamma(\delta + 1)|p_m|}\right) \\
&\leqslant \sum_{k \geqslant 0} \sum_{d \geqslant 1} M_{k,d} 2^k \exp\left(-\frac{(1+d) \log^{1+\epsilon}\left(2^k(1+d)\right)}{10\theta\gamma(\delta+1)}\right) \\
&\leqslant A R_2' R_2''
\end{aligned}
$$

where we have set

$$
R_2' = \sum_{k \geqslant 0} \exp\left(k \log 2 - \frac{(k \log 2)^{1+\epsilon}}{5\theta\gamma(\delta+1)}\right) < \infty
$$

and

$$
R_2'' = \sum_{d \geqslant 1} \exp\left(B \log(1+d) - \frac{(1+d)\log^{1+\epsilon}(1+d)}{10\theta\gamma(\delta+1)}\right) < \infty .
$$

We now have to prove theorem 4. For an arbitrary $m \in \mathcal{M}$, we begin the proof by expanding the Kullback-Leibler divergence of $(\tilde{s}, \tilde{\sigma})$,

$$
\begin{aligned}
\mathcal{K}(P_{s,\sigma}, P_{\tilde{s},\tilde{\sigma}}) &= \frac{1}{2} \sum_{i=1}^n \frac{(s_i - \tilde{s}_i)^2}{\tilde{\sigma}_i} + \phi\left(\frac{\tilde{\sigma}_i}{\sigma_i}\right) \\
&= \mathcal{K}(P_{s,\sigma}, P_{\hat{s}_m, \hat{\sigma}_m}) + [\mathcal{L}(\hat{s}_m, \hat{\sigma}_m) - \mathcal{K}(P_{s,\sigma}, P_{\hat{s}_m, \hat{\sigma}_m})] \\
&\quad + [\mathcal{L}(\tilde{s}, \tilde{\sigma}) - \mathcal{L}(\hat{s}_m, \hat{\sigma}_m)] + [\mathcal{K}(P_{s,\sigma}, P_{\tilde{s},\tilde{\sigma}}) - \mathcal{L}(\tilde{s}, \tilde{\sigma})] .
\end{aligned}
$$

By the definition (2.3) of $\hat{m}$, the inequality

$$
\mathcal{L}(\tilde{s}, \tilde{\sigma}) - \mathcal{L}(\hat{s}_m, \hat{\sigma}_m) \leqslant \text{pen}(m) - \text{pen}(\hat{m}) \tag{4.7}
$$

is true for any $m \in \mathcal{M}$. The difference between the divergence and the likelihood can be expressed as

$$
\begin{aligned}
&\mathcal{K}(P_{s,\sigma}, P_{\hat{s}_m, \hat{\sigma}_m}) - \mathcal{L}(\hat{s}_m, \hat{\sigma}_m) \\
&= \frac{1}{2} \sum_{I \in p_m} \sum_{i \in I} \left(\frac{\sigma_i}{\hat{\sigma}_{m,I}} - 1\right)\left(1 - \varepsilon_i^{[1]^2}\right) \tag{4.8} \\
&\qquad - \frac{2(s_i - \hat{s}_{m,i})\sqrt{\sigma_i}\varepsilon_i^{[1]}}{\hat{\sigma}_{m,I}} - \frac{1}{2} \sum_{i=1}^n \left(\varepsilon_i^{[1]^2} + \log \sigma_i\right) .
\end{aligned}
$$

Using (4.7) and (4.8), for any $\alpha \in (0,1)$, we can write

$$(1-\alpha)\mathcal{K}(P_{s,\sigma}, P_{\tilde{s},\tilde{\sigma}}) \tag{4.9}$$
$$\leqslant \mathcal{K}(P_{s,\sigma}, P_{\hat{s}_m, \hat{\sigma}_m}) + \text{pen}(m) + G(m)$$
$$+ W_1(\hat{m}) + W_2(\hat{m}) + Z(\hat{m}) - \text{pen}(\hat{m})$$

where, for any $m \in \mathcal{M}$,

$$W_1(m) = \sum_{I \in p_m} \frac{1}{\hat{\sigma}_{m,I}} \left\| \pi_m \Gamma_\sigma^{1/2} \varepsilon^{[1]} \right\|_I^2,$$

$$W_2(m) = \sum_{I \in p_m} \frac{1}{\hat{\sigma}_{m,I}} \left( \left\langle s_m - s, \Gamma_\sigma^{1/2} \varepsilon^{[1]} \right\rangle_I - \frac{\alpha}{2} \|s_m - s\|_I^2 \right),$$

$$Z(m) = \frac{1}{2} \sum_{I \in p_m} \sum_{i \in I} \left( \left( \frac{\sigma_i}{\hat{\sigma}_{m,I}} - 1 \right) \left( 1 - \varepsilon_i^{[1]^2} \right) - \alpha \phi \left( \frac{\hat{\sigma}_{m,I}}{\sigma_i} \right) \right)$$

and

$$G(m) = \sum_{I \in p_m} \left( \frac{1}{\hat{\sigma}_{m,I}} \left\langle s - \hat{s}_m, \Gamma_\sigma^{1/2} \varepsilon^{[1]} \right\rangle_I - \frac{1}{2} \sum_{i \in I} \left( \frac{\sigma_i}{\hat{\sigma}_{m,I}} - 1 \right) \left( 1 - \varepsilon_i^{[1]^2} \right) \right).$$

We split the proof of theorem 4 in several lemmas.

**Lemma 5.** *For any $m \in \mathcal{M}$, we have*

$$\mathbb{E}[G(m)] \leqslant 0.$$

*Proof.* Let us compute this expectation to obtain the inequality. By independence between $\varepsilon^{[1]}$ and $\varepsilon^{[2]}$, we get

$$\mathbb{E}\left[G(m) \big| \varepsilon^{[2]}\right] = \sum_{I \in p_m} -\frac{1}{\hat{\sigma}_{m,I}} \mathbb{E}\left[\left\langle \pi_m \Gamma_\sigma^{1/2} \varepsilon^{[1]}, \Gamma_\sigma^{1/2} \varepsilon^{[1]} \right\rangle_I\right]$$
$$= -\sum_{I \in p_m} \frac{1}{\hat{\sigma}_{m,I}} \mathbb{E}\left[\left\| \pi_m \Gamma_\sigma^{1/2} \varepsilon^{[1]} \right\|_I^2\right].$$

It leads to $\mathbb{E}[G(m)] = \mathbb{E}\left[\mathbb{E}\left[G(m) \big| \varepsilon^{[2]}\right]\right] \leqslant 0$. $\square$

In order to control $Z(m)$, we split it in two terms that we study separately,

$$Z(m) = Z_+(m) + Z_-(m)$$

where

$$Z_+(m) = \frac{1}{2} \sum_{I \in p_m} \sum_{i \in I} \left( \left( \frac{\sigma_i}{\hat{\sigma}_{m,I}} - 1 \right)_+ \left( 1 - \varepsilon_i^{[1]^2} \right) - \alpha \phi \left( \frac{\hat{\sigma}_{m,I}}{\sigma_i} \right) \mathbb{1}_{\hat{\sigma}_{m,I} \leqslant \sigma_i} \right)$$

and

$$Z_-(m) = \frac{1}{2} \sum_{I \in p_m} \sum_{i \in I} \left( \left( \frac{\sigma_i}{\hat{\sigma}_{m,I}} - 1 \right)_- \left( \varepsilon_i^{[1]^2} - 1 \right) - \alpha \phi \left( \frac{\hat{\sigma}_{m,I}}{\sigma_i} \right) \mathbb{1}_{\hat{\sigma}_{m,I} > \sigma_i} \right).$$

**Lemma 6.** *Let $m \in \mathcal{M}$ and $x$ be a positive number. Under the hypothesis $(\mathbf{H}_\theta)$, we get*

$$\mathbb{E}\left[(Z_+(m) - x)_+\right] \leqslant \frac{\gamma\theta|p_m|}{\alpha} \exp\left(-\frac{n - (d_m + 3)|p_m|}{2|p_m|} \log\left(1 + \frac{2\alpha|p_m|x}{\gamma n}\right)\right) .$$

*Proof.* We begin by setting, for all $1 \leqslant i \leqslant n$,

$$T_i(m) = \frac{(\sigma_i/\hat{\sigma}_{m,i} - 1)_+}{\left(\sum_{j=1}^n (\sigma_j/\hat{\sigma}_{m,j} - 1)_+^2\right)^{1/2}}$$

and we denote by

$$S(m) = \sum_{i=1}^n T_i(m) \left(1 - \varepsilon_i^{[1]\,2}\right) .$$

We lower bound the function $\phi$ by the remark

$$\forall a \in (0,1),\ \forall u \in [a,1],\ \left(\frac{1}{u} - 1\right)^2 \leqslant \frac{2}{a} \phi(u) .$$

Thus, we obtain

$$\sum_{i=1}^n \left(\frac{\sigma_i}{\hat{\sigma}_{m,i}} - 1\right)_+^2 \leqslant 2 \left(\max_{i \leqslant n} \frac{\sigma_i}{\hat{\sigma}_{m,i}}\right) \sum_{j=1}^n \phi\left(\frac{\hat{\sigma}_{m,j}}{\sigma_j}\right) \mathbb{1}_{\hat{\sigma}_{m,j} \leqslant \sigma_j} = 2M(m)$$

and we use this inequality to get

$$\begin{aligned}
Z_+(m) &= \frac{1}{2} \left(\sum_{i=1}^n \left(\frac{\sigma_i}{\hat{\sigma}_{m,i}} - 1\right)_+^2\right)^{1/2} S(m) - \frac{\alpha}{2} \sum_{i=1}^n \phi\left(\frac{\hat{\sigma}_{m,i}}{\sigma_i}\right) \mathbb{1}_{\hat{\sigma}_{m,i} \leqslant \sigma_i} \\
&\leqslant \sqrt{\frac{M(m)}{2}} S(m)_+ - \frac{\alpha}{2} \sum_{i=1}^n \phi\left(\frac{\hat{\sigma}_{m,i}}{\sigma_i}\right) \mathbb{1}_{\hat{\sigma}_{m,i} \leqslant \sigma_i} \\
&\leqslant \frac{1}{4\alpha} \left(\max_{i \leqslant n} \frac{\sigma_i}{\hat{\sigma}_{m,i}}\right) S(m)_+^2 .
\end{aligned}$$

To control $S(m)$, we use the inequality (4.2) in [15], conditionally to $\varepsilon^{[2]}$. Let $u > 0$,

$$\begin{aligned}
\mathbb{P}\left(\left(\max_{i \leqslant n} \frac{\sigma_i}{\hat{\sigma}_{m,i}}\right) S(m)_+^2 \geqslant u\right) &= \mathbb{E}\left[\mathbb{P}\left(S(m) \geqslant \sqrt{u/\max_{i \leqslant n} \frac{\sigma_i}{\hat{\sigma}_{m,i}}}\,\bigg|\,\varepsilon^{[2]}\right)\right] \\
&\leqslant \mathbb{E}\left[\exp\left(-\frac{u}{4} \min_{i \leqslant n} \frac{\hat{\sigma}_{m,i}}{\sigma_i}\right)\right] .
\end{aligned}$$

By the remark (4.1), we can upper bound it by

$$\mathbb{P}\left(\left(\max_{i \leqslant n} \frac{\sigma_i}{\hat{\sigma}_{m,i}}\right) S(m)_+^2 \geqslant u\right) \leqslant \mathbb{E}\left[\exp\left(-\frac{u}{4\gamma} \min_{I \in p_m} X_I\right)\right]$$

where the $X_I$'s are *i.i.d.* random variables with a $\chi^2\left(|I| - d_m - 1\right)/|I|$ distribution.

For any $\lambda > 0$, we know that the Laplace transform of $X_I$ is given by

$$\mathbb{E}\left[e^{-\lambda X_I}\right] = \left(1 + \frac{2\lambda}{|I|}\right)^{-(|I|-d_m-1)/2}. \tag{4.10}$$

Let $t > 0$, the following expectation is dominated by

$$\begin{aligned}
\mathbb{E}\left[\left(Z_+(m) - \frac{\gamma}{2\alpha}t\right)_+\right] &= \int_0^\infty \mathbb{P}\left(Z_+(m) \geqslant \frac{\gamma t}{2\alpha} + u\right) du \\
&\leqslant \int_0^\infty \mathbb{E}\left[\exp\left(-\left(\frac{\alpha u}{\gamma} + \frac{t}{2}\right) \min_{I \in p_m} X_I\right)\right] du \\
&\leqslant \int_0^\infty \mathbb{E}\left[\max_{I \in p_m} \exp\left(-\left(\frac{\alpha u}{\gamma} + \frac{t}{2}\right) X_I\right)\right] du.
\end{aligned}$$

Using $(\mathbf{H}_\theta)$ and (4.10), we roughly upper bound the maximum by the sum of the Laplace transforms and we get

$$\begin{aligned}
\mathbb{E}&\left[\left(Z_+(m) - \frac{\gamma}{2\alpha}t\right)_+\right] \\
&\leqslant \sum_{I \in p_m} \frac{\gamma|I|}{\alpha(|I| - d_m - 3)}\left(1 + \frac{t}{|I|}\right)^{-(|I|-d_m-3)/2} \\
&\leqslant \frac{\gamma\theta|p_m|}{\alpha} \exp\left(-\frac{n - (d_m+3)|p_m|}{2|p_m|}\log\left(1 + \frac{t|p_m|}{n}\right)\right).
\end{aligned}$$

Take $t = 2\alpha x/\gamma$ to conclude. $\square$

**Lemma 7.** *Let $m \in \mathcal{M}$ and $x$ be a positive number, then*

$$\mathbb{E}\left[(Z_-(m) - (2\alpha+1)x)_+\right] \leqslant \frac{2\alpha+1}{\alpha}e^{-\alpha x}.$$

*Proof.* Note that for all $u > 1$, we have

$$2\phi(u) \geqslant \left(\frac{1}{u} - 1\right)^2.$$

Let $t > 0$, we handle $Z_-(m)$ conditionally to $\varepsilon^{[2]}$ and, using the previous lower bound on $\phi$, we obtain

$$\mathbb{P}\left(Z_-(m) \geqslant \frac{2\alpha+1}{2\alpha}t \,\bigg|\, \varepsilon^{[2]}\right)$$

$$\leqslant \mathbb{P}\left(\frac{1}{2}\sum_{i=1}^n \left(\frac{\sigma_i}{\hat{\sigma}_{m,i}} - 1\right)_- \left(\varepsilon_i^{[1]2} - 1\right) \geqslant \frac{2\alpha+1}{2\alpha}t + \frac{\alpha}{4}\sum_{i=1}^n \left(\frac{\sigma_i}{\hat{\sigma}_{m,i}} - 1\right)_-^2 \,\bigg|\, \varepsilon^{[2]}\right)$$

$$\leqslant \mathbb{P}\left(\frac{1}{2}\sum_{i=1}^{n}\left(\frac{\sigma_i}{\hat\sigma_{m,i}}-1\right)_{-}\left(\varepsilon_i^{[1]2}-1\right) \geqslant t + \sqrt{\frac{t}{2}\sum_{i=1}^{n}\left(\frac{\sigma_i}{\hat\sigma_{m,i}}-1\right)_{-}^{2}}\,\bigg|\varepsilon^{[2]}\right).$$

Let us note that

$$\max_{i\leqslant n}\left(\frac{\sigma_i}{\hat\sigma_{m,i}}-1\right)_{-} \leqslant 1\,,$$

thus, we can apply the inequality (4.1) from [15] to get

$$\mathbb{P}\left(Z_{-}(m)\geqslant \frac{2\alpha+1}{2\alpha}t\right) \leqslant \exp(-t/2)\,.$$

This inequality leads us to

$$\mathbb{E}\left[\left(Z_{-}(m)-\frac{2\alpha+1}{\alpha}t\right)_{+}\right] \leqslant \int_{(2\alpha+1)t/\alpha}^{+\infty}\mathbb{P}(Z_{-}(m)\geqslant u)du$$
$$\leqslant \frac{2\alpha+1}{\alpha}e^{-t}\,.$$

Take $t=\alpha x$ to get the announced result. □

It remains to control $W_1(m)$ and $W_2(m)$. For the first one, we now prove a Rosenthal-type inequality.

**Lemma 8.** *Consider any $m\in\mathcal{M}$. Under the hypothesis $(\mathbf{H}_\theta)$, for any $x>0$, we have*

$$\mathbb{E}[(W_1(m)-\gamma\theta D_m-x)_{+}]$$
$$\leqslant C\theta^2\gamma\sqrt{|p_m|}d_m\left(\frac{2C\theta^2\gamma\sqrt{|p_m|}d_m\log(1+d_m)}{x}\right)^{\lfloor 2\log(1+d_m)\rfloor}$$

*where $\lfloor\cdot\rfloor$ is the integral part and $C$ is a positive constant that could be taken equal to*

$$C=\frac{12\sqrt{2e}}{\sqrt{e}-1}\approx 43.131\,.$$

*Proof.* Using the lemma 10 and the remark (4.1), we dominate $W_1(m)$,

$$W_1(m)\leqslant W_1'(m)=\gamma\sum_{I\in p_m}\frac{|I|d_m}{|I|-d_m-1}F_I = \frac{\gamma n d_m}{n-|p_m|(1+d_m)}\sum_{I\in p_m}F_I$$

where the $F_I$'s are *i.i.d.* Fisher random variables of parameters $(d_m, n/|p_m|-d_m-1)$. We denote by $F_m$ the distribution of the $F_I$'s and we have

$$\frac{\gamma}{2}D_m \leqslant \gamma|p_m|d_m \leqslant \mathbb{E}[W_1'(m)]\leqslant \gamma\theta|p_m|d_m \leqslant \gamma\theta D_m\,.$$

Take $x>0$ and an integer $q>1$, then

$$\mathbb{E}\left[(W_1'(m)-\mathbb{E}[W_1'(m)]-x)_{+}\right] \leqslant \frac{\mathbb{E}\left[(W_1'(m)-\mathbb{E}[W_1'(m)])_{+}^{q}\right]}{(q-1)x^{q-1}}\,. \qquad (4.11)$$

We set $V = W'_1(m) - \mathbb{E}[W'_1(m)]$. It is the sum of the independent centered random variables

$$X_I = \frac{\gamma n d_m}{n - |p_m|(1+d_m)}(F_I - \mathbb{E}[F_I]), \ I \in p_m \ .$$

To dominate $\mathbb{E}\left[V_+^q\right]$, we use the theorem 9 in [11]. Let us compute

$$\sum_{I \in p_m} \mathbb{E}[X_I^2] = \frac{2\gamma^2 n^2 d_m (n - 3|p_m|)|p_m|}{(n - |p_m|(d_m + 3))^2 (n - |p_m|(d_m + 5))} \leqslant 2\gamma^2 \theta^3 |p_m| d_m$$

and so,

$$\mathbb{E}\left[V_+^q\right]^{1/q} \leqslant \sqrt{12\kappa' \gamma^2 \theta^3 |p_m| d_m q} + q\kappa' \sqrt{2} \mathbb{E}\left[\max_{I \in p_m} |X_I|^q\right]^{1/q}$$

where $\kappa' = \frac{\sqrt{e}}{2(\sqrt{e}-1)}$.

We consider $q = 1 + \lfloor 2\log(1 + d_m) \rfloor$ where $\lfloor \cdot \rfloor$ is the integral part. For this choice, $q \leqslant 1 + d_m$ and it implies

$$2|p_m|q < n - |p_m|(1 + d_m) \ .$$

The hypothesis $(\mathbf{H}_\theta)$ allows us to make a such choice. We roughly upper bound the maximum by the sum and we use $(\mathbf{H}_\theta)$ to get

$$\begin{aligned}
\mathbb{E}\left[\max_{I \in p_m}|X_I|^q\right] &\leqslant (\gamma \theta d_m)^q \mathbb{E}\left[\max_{I \in p_m}|F_I - \mathbb{E}[F_I]|^q\right] \\
&\leqslant (\gamma \theta d_m)^q 2^{q-1} \left(\mathbb{E}[F_m]^q + |p_m|\mathbb{E}[F_m^q]\right) \\
&\leqslant \frac{(2\gamma \theta^2 d_m)^q}{2} + \frac{|p_m|}{2}\left(\frac{(2\gamma \theta d_m)(1 + 2(q-1)/d_m)}{1 - 2|p_m|q/(n - |p_m|(1+d_m))}\right)^q \\
&\leqslant (6\gamma \theta^2 d_m)^q |p_m| \ .
\end{aligned}$$

Thus, it gives

$$\begin{aligned}
\mathbb{E}\left[V_+^q\right]^{1/q} &\leqslant \gamma \theta^2 \left(\sqrt{12\kappa' |p_m| d_m q} + 6\kappa' \sqrt{2}|p_m|^{1/q} d_m q\right) \\
&\leqslant 6\kappa' \sqrt{2} \gamma \theta^2 \left(\sqrt{|p_m| d_m q} + |p_m|^{1/q} d_m q\right) \\
&\leqslant 12\kappa' \sqrt{2} \gamma \theta^2 \sqrt{|p_m| d_m} (1 + \lfloor 2\log(1+d_m)\rfloor) \ .
\end{aligned}$$

Injecting this inequality in (4.11) leads to

$$\mathbb{E}\left[(W'_1(m) - \mathbb{E}[W'_1(m)] - x)_+\right]$$
$$\leqslant C\gamma \theta^2 \sqrt{|p_m|} d_m \left(\frac{C\gamma \theta^2 \sqrt{|p_m|} d_m (1 + 2\log(1+d_m))}{2x}\right)^{\lfloor 2\log(1+d_m) \rfloor} .$$

□

**Lemma 9.** *Consider any $m \in \mathcal{M}$ and let $x$ be a positive number. Under the hypothesis $(H_\theta)$, we have*

$$\mathbb{E}\left[(W_2(m) - x)_+\right] \leqslant \frac{\gamma\theta|p_m|}{\alpha} \exp\left(-\frac{n - (d_m + 3)|p_m|}{2|p_m|} \log\left(1 + \frac{2\alpha|p_m|x}{\gamma n}\right)\right).$$

*Proof.* Let us define

$$A(m) = \sum_{I \in p_m} \frac{\|s - s_m\|_I^2}{\hat{\sigma}_{m,I}}.$$

The distribution of $W_2(m)$ conditionally to $\varepsilon^{[2]}$ is Gaussian with mean equal to $-\alpha A(m)/2$ and variance factor

$$\sum_{I \in p_m} \frac{\left\|\Gamma_\sigma^{1/2}(s - s_m)\right\|_I^2}{\hat{\sigma}_{m,I}^2}.$$

If $\zeta$ is a standard Gaussian random variable, it is well known that, for any $\lambda > 0$,

$$\mathbb{P}(\zeta \geqslant \sqrt{2\lambda}) \leqslant e^{-\lambda}. \tag{4.12}$$

We apply the Gaussian inequality (4.12) to $W_2(m)$ conditionally to $\varepsilon^{[2]}$,

$$\forall t > 0, \ \mathbb{P}\left(W_2(m) + \frac{\alpha}{2}A(m) \geqslant \sqrt{2t \sum_{I \in p_m} \frac{\left\|\Gamma_\sigma^{1/2}(s - s_m)\right\|_I^2}{\hat{\sigma}_{m,I}^2}} \ \Big|\ \varepsilon^{[2]}\right) \leqslant e^{-t}.$$

It leads to

$$\mathbb{P}\left(W_2(m) + \frac{\alpha}{2}A(m) \geqslant \sqrt{2tA(m) \max_{i \leqslant n} \frac{\sigma_i}{\hat{\sigma}_{m,i}}} \ \Big|\ \varepsilon^{[2]}\right) \leqslant e^{-t}$$

and thus, by the remark (4.1),

$$\mathbb{P}\left(W_2(m) \geqslant \frac{\gamma t}{\alpha} \max_{I \in p_m} X_I^{-1} \ \Big|\ \varepsilon^{[2]}\right) \leqslant \mathbb{P}\left(W_2(m) \geqslant \frac{t}{\alpha} \max_{i \leqslant n} \frac{\sigma_i}{\hat{\sigma}_{m,i}} \ \Big|\ \varepsilon^{[2]}\right) \leqslant e^{-t}$$

where the $X_I$'s are *i.i.d.* random variables with a $\chi^2(|I| - d_m - 1)/|I|$ distribution. Finally, we integrate following $\varepsilon^{[2]}$ and we get

$$\mathbb{P}(W_2(m) \geqslant t) \leqslant \mathbb{E}\left[\max_{I \in p_m} \exp\left(-\frac{\alpha t}{\gamma} X_I\right)\right].$$

We finish as we did for $Z_+(m)$,

$$\mathbb{E}\left[\left(W_2(m) - \frac{\gamma}{2\alpha}t\right)_+\right]$$

$$\leqslant \int_0^{+\infty} \mathbb{E}\left[\max_{I \in p_m} \exp\left(-\left(\frac{\alpha u}{\gamma} + \frac{t}{2}\right)X_I\right)\right] du$$

$$\leqslant \frac{\gamma\theta}{\alpha} \sum_{I \in p_m} \left(1 + \frac{t}{|I|}\right)^{-(|I| - d_m - 3)/2}$$

$$\leqslant \frac{\gamma\theta|p_m|}{\alpha} \exp\left(-\frac{n - (d_m + 3)|p_m|}{2|p_m|} \log\left(1 + \frac{t|p_m|}{n}\right)\right).$$

$\square$

In order to end the proof of theorem 4, we need to put together the results of the previous lemmas. Because $\gamma \geqslant 1$, for any $x > 0$, we can write

$$e^{-\alpha x} \leqslant \exp\left(-\frac{n}{2|p_m|}\log\left(1+\frac{2\alpha|p_m|x}{\gamma n}\right)\right).$$

We now come back to (4.9) and we apply the preceding results to each model. Let $m \in \mathcal{M}$, we take

$$x = \frac{x_m}{2(2+\alpha)}$$

and, recalling (4.6), we get the following inequalities

$$(1-\alpha)\mathbb{E}[\mathcal{K}(P_{s,\sigma}, P_{\tilde{s},\tilde{\sigma}})]$$
$$\leqslant \mathbb{E}[\mathcal{K}(P_{s,\sigma}, P_{\hat{s}_m,\hat{\sigma}_m})] + \mathrm{pen}(m) + \mathbb{E}\left[\left(W_1(\hat{m}) - \gamma\theta D_{\hat{m}} - \frac{x_{\hat{m}}}{2(2+\alpha)}\right)_+\right]$$
$$+ \mathbb{E}\left[\left(W_2(\hat{m}) - \frac{x_{\hat{m}}}{2(2+\alpha)}\right)_+\right] + \mathbb{E}\left[\left(Z_+(\hat{m}) - \frac{x_{\hat{m}}}{2(2+\alpha)}\right)_+\right]$$
$$+ \mathbb{E}\left[\left(Z_-(\hat{m}) - (1+2\alpha)\frac{x_{\hat{m}}}{2(2+\alpha)}\right)_+\right]$$
$$\leqslant \mathbb{E}[\mathcal{K}(m)] + \mathrm{pen}(m) + R_1(\mathcal{M}) + R_2(\mathcal{M}) \qquad (4.13)$$

where $R_1(\mathcal{M})$ and $R_2(\mathcal{M})$ are the sums defined in the theorem 4. As the choice of $m$ is arbitrary, we can take the infimum among $m \in \mathcal{M}$ in the right part of (4.13).

### 4.3. Proof of the proposition 3

For the collection $\mathcal{F}^{PC}$, we have $A = 1$ and $B = 0$ in (2.4). Let $m \in \mathcal{M}$, we denote by $\bar{\sigma}_m \in \Sigma_m$ the quantity

$$\bar{\sigma}_m = \sum_{I \in p_m} \bar{\sigma}_{m,I}\mathbb{1}_I \text{ with } \forall I \in p_m, \ \bar{\sigma}_{m,I} = \frac{1}{|I|}\sum_{i \in I}\sigma_i.$$

The theorem 2 gives us

$$\mathbb{E}\left[\mathcal{K}_n(P_{s,\sigma}, P_{\tilde{s},\tilde{\sigma}})\right]$$
$$\leqslant \frac{C}{n}\inf_{m \in \mathcal{M}}\left\{\mathcal{K}(P_{s,\sigma}, P_{s_m,\sigma_m}) + D_m\log^{1+\epsilon}D_m\right\} + \frac{R}{n}$$
$$\leqslant \frac{C}{n}\inf_{m \in \mathcal{M}}\left\{\mathcal{K}(P_{s,\sigma}, P_{s_m,\bar{\sigma}_m}) + D_m\log^{1+\epsilon}D_m\right\} + \frac{R}{n}$$
$$\leqslant C\inf_{m \in \mathcal{M}}\left\{\frac{\|s-s_m\|_2^2}{2n\sigma_*} + \frac{\|\sigma-\bar{\sigma}_m\|_2^2}{2n\sigma_*^2} + D_m\log^{1+\epsilon}D_m\right\} + \frac{R}{n}$$

because, for any $x > 0$, $\phi(x) \leqslant (x - 1/x)^2$.

Assuming $(s_r, \sigma_r) \in \mathcal{H}_{\alpha_1}(L_1) \times \mathcal{H}_{\alpha_2}(L_2)$, we know (see [13]) that

$$\|s - s_m\|_2^2 \leqslant n L_1^2 (|p_m| d_m)^{-2\alpha_1}$$

and

$$\|\sigma - \bar{\sigma}_m\|_2^2 \leqslant n L_2^2 |p_m|^{-2\alpha_2} .$$

Thus, we obtain

$$\mathbb{E}\left[\mathcal{K}_n(P_{s,\sigma}, P_{\tilde{s},\tilde{\sigma}})\right]$$
$$\leqslant C \inf_{m \in \mathcal{M}} \left\{ \frac{L_1^2}{2\sigma_*}(|p_m|d_m)^{-2\alpha_1} + \frac{L_2^2}{2\sigma_*^2}|p_m|^{-2\alpha_2} + \frac{\log^{1+\epsilon} n}{n} D_m \right\} + \frac{R}{n} .$$

If $\alpha_1 < \alpha_2$, we can take

$$|p_m| d_m = \left\lfloor \left( \frac{L_1^2 n}{2\sigma_* \log^{1+\epsilon} n} \right)^{1/(1+2\alpha_1)} \right\rfloor$$

and

$$|p_m| = \left\lfloor \left( \frac{L_2^2 n}{2\sigma_*^2 \log^{1+\epsilon} n} \right)^{1/(1+2\alpha_2)} \right\rfloor .$$

For $\alpha_1 \geqslant \alpha_2$, this choice is not allowed because it would imply $d_m = 0$. So, in this case, we take

$$d_m = 1 \text{ and } |p_m| = \left\lfloor \left( \frac{(L_1^2 \sigma_* + L_2^2) n}{2\sigma_*^2 \log^{1+\epsilon} n} \right)^{1/(1+2\alpha_2)} \right\rfloor .$$

In the two situation, we obtain the announced result.

## 5. Technical results

This section is devoted to some useful technical results. Some notations previously introduced can have a different meaning here.

**Lemma 10.** *Let $\Sigma$ be a positive symmetric $n \times n$-matrix and $\sigma_1, \ldots, \sigma_n > 0$ be its eigenvalues. Let $P$ be an orthogonal projection of rank $D \geqslant 1$. If we denote $M = P\Sigma P$, then $M$ is a non-negative symmetric matrix of rank $D$ and, if $\tau_1, \ldots, \tau_D$ are its positive eigenvalues, we have*

$$\min_{1 \leqslant i \leqslant n} \sigma_i \leqslant \min_{1 \leqslant i \leqslant D} \tau_i \quad \text{and} \quad \max_{1 \leqslant i \leqslant D} \tau_i \leqslant \max_{1 \leqslant i \leqslant n} \sigma_i .$$

*Proof.* We denote by $\Sigma^{1/2}$ the symmetric square root of $\Sigma$. By a classical result, $M$ has the same rank, equal to $D$, than $P\Sigma^{1/2}$. On a first side, we have

$$\max_{1 \leqslant i \leqslant D} \tau_i = \sup_{\substack{x \in \mathbb{R}^n \\ x \neq 0}} \frac{\langle P\Sigma P x, \, x \rangle}{\|x\|^2}$$

$$= \sup_{\substack{(x_1, x_2) \in \ker(P) \times \operatorname{im}(P) \\ (x_1, x_2) \neq (0,0)}} \frac{\langle P\Sigma x_2, \, x_2 \rangle}{\|x_1\|^2 + \|x_2\|^2}$$

$$\leqslant \sup_{\substack{x_2 \in \operatorname{im}(P) \\ x_2 \neq 0}} \frac{\langle \Sigma x_2, \, x_2 \rangle}{\|x_2\|^2} \leqslant \max_{1 \leqslant i \leqslant n} \sigma_i \; .$$

On the other side, we can write

$$\min_{1 \leqslant i \leqslant D} \tau_i = \min_{\substack{V \subset \mathbb{R}^n \\ \dim(V) = n - D + 1}} \max_{\substack{x \in V \\ x \neq 0}} \frac{\langle M x, \, x \rangle}{\|x\|^2}$$

$$= \min_{\substack{V \subset \mathbb{R}^n \\ \dim(V) = n - D + 1}} \max_{\substack{x \in V \\ x \neq 0}} \frac{\|\Sigma^{1/2} P x\|^2}{\|x\|^2}$$

$$\geqslant \min_{\substack{V \subset \mathbb{R}^n \\ \dim(V) = n - D + 1}} \max_{\substack{x \in V \cap \operatorname{im}(P) \\ x \neq 0}} \frac{\|\Sigma^{1/2} x\|^2}{\|x\|^2}$$

$$\geqslant \min_{\substack{V' \subset \mathbb{R}^n \\ \dim(V') \geqslant 1}} \max_{\substack{x \in V' \\ x \neq 0}} \frac{\|\Sigma^{1/2} x\|^2}{\|x\|^2} = \min_{1 \leqslant i \leqslant n} \sigma_i \; .$$

$\square$

**Lemma 11.** *Let $\varepsilon$ be a standard Gaussian vector in $\mathbb{R}^n$, $a = (a_1, \ldots, a_n)' \in \mathbb{R}^n$ and $b_1, \ldots, b_n > 0$. We denote by $b^*$ (resp. $b_*$) the maximum (resp. minimum) of the $b_i$'s. If $n > 2$ and $Z = \sum_{i=1}^n (a_i + \sqrt{b_i} \varepsilon_i)^2$, then*

$$\mathbb{E}\left[\frac{1}{Z}\right] \leqslant \frac{1}{\mathbb{E}[Z]} \left( 1 + \frac{2\kappa (b^*/b_*)^2}{n-2} \right)$$

*where $\kappa > 1$ is a constant that can be taken equal to $1 + 2e^{-1} \approx 1.736$.*

*Proof.* We recall that $\mathbb{E}[Z] = \sum_{i=0}^n (a_i^2 + b_i)$ and, for any $\lambda > 0$, the Laplace transform of $(a_i + \sqrt{b_i} \varepsilon_i)^2$ is

$$\mathbb{E}\left[ \exp\left( -\lambda (a_i + \sqrt{b_i} \varepsilon_i)^2 \right) \right] = \exp\left( -\frac{\lambda a_i^2}{1 + 2\lambda b_i} - \frac{1}{2} \log(1 + 2\lambda b_i) \right) \; .$$

Thus, the Laplace transform of $Z$ is equal to

$$\psi(\lambda) = \mathbb{E}\left[ e^{-\lambda Z} \right]$$

$$= \exp\left( -\sum_{i=1}^n \frac{\lambda a_i^2}{1 + 2\lambda b_i} - \frac{1}{2} \sum_{i=1}^n \log(1 + 2\lambda b_i) \right)$$

$$= e^{-\lambda \mathbb{E}[Z]} \times \exp\left( \sum_{i=0}^n \frac{2\lambda^2 a_i^2 b_i}{1 + 2\lambda b_i} - \frac{1}{2} \sum_{i=1}^n r(2\lambda b_i) \right)$$

where $r(x) = \log(1+x) - x$ for all $x > 0$. To compute the expectation of the inverse of $Z$, we integrate $\psi$ by parts,

$$\begin{aligned}
\mathbb{E}\left[\frac{1}{Z}\right] &= \int_0^\infty \psi(\lambda) d\lambda \\
&= \int_0^\infty e^{-\lambda \mathbb{E}[Z]} \times \exp\left(\sum_{i=0}^n \frac{2\lambda^2 a_i^2 b_i}{1+2\lambda b_i} - \frac{1}{2}\sum_{i=1}^n r(2\lambda b_i)\right) d\lambda \\
&= \frac{1}{\mathbb{E}[Z]} + \frac{1}{\mathbb{E}[Z]} \int_0^\infty f_{a,b}(\lambda) \psi(\lambda) d\lambda
\end{aligned}$$

where

$$f_{a,b}(\lambda) = \sum_{i=0}^n \frac{2\lambda b_i^2}{1+2\lambda b_i} + \frac{4\lambda a_i^2 b_i (1+\lambda b_i)}{(1+2\lambda b_i)^2} \ .$$

We now upper bound the integral,

$$\begin{aligned}
\mathbb{E}\left[\frac{\mathbb{E}[Z]}{Z} - 1\right] &= \int_0^\infty f_{a,b}(\lambda) \frac{\exp\left(-\sum_{i=1}^n \lambda a_i^2/(1+2\lambda b_i)\right)}{\prod_{i=1}^n \sqrt{1+2\lambda b_i}} d\lambda \\
&\leqslant \int_0^\infty \frac{2n\lambda b^{*2}}{(1+2\lambda b_*)^{1+n/2}} d\lambda \\
&\quad + \int_0^\infty \frac{4b^*(1+\lambda b^*)}{(1+2\lambda b_*)^{1+n/2}} \times g_{a,b}(\lambda) e^{-g_{a,b}(\lambda)} d\lambda
\end{aligned}$$

where we have set

$$g_{a,b}(\lambda) = \sum_{i=1}^n \frac{\lambda a_i^2}{1+2\lambda b_i} \ .$$

For any $t > 0$, $te^{-t} \leqslant e^{-1}$. Because $g_{a,b}$ is a positive function and $n > 2$, we obtain

$$\begin{aligned}
\mathbb{E}\left[\frac{\mathbb{E}[Z]}{Z} - 1\right] &\leqslant \int_0^\infty \frac{2n\lambda b^{*2}}{(1+2\lambda b_*)^{1+n/2}} d\lambda + \int_0^\infty \frac{4b^*(1+\lambda b^*)}{e(1+2\lambda b_*)^{1+n/2}} d\lambda \\
&\leqslant \frac{2(b^*/b_*)^2}{n-2} + \frac{4(b^*/b_*)(n-2+b^*/b_*)}{en(n-2)} \\
&\leqslant \frac{2(b^*/b_*)^2}{n-2}\left(1 + \frac{2(n-1)}{en}\right) \\
&\leqslant 2(1+2e^{-1})\frac{(b^*/b_*)^2}{n-2} \ .
\end{aligned}$$

$\square$